\documentclass{article}

\usepackage[english]{babel}\usepackage[latin1]{inputenc}\usepackage[T1]{fontenc}
\usepackage{graphicx,amsfonts,amssymb,amsmath,amsthm,mathrsfs,enumitem,color}
\usepackage[top=2cm,bottom=3cm,left=2cm,right=2cm]{geometry}
\usepackage{sectsty}\sectionfont{\normalsize}\subsectionfont{\normalsize}

\newtheorem{theorem}{Theorem}[section]

\newtheorem{lemma}{Lemma}[section]
\newtheorem{corollary}{Corollary}[section]
\theoremstyle{definition}
\newtheorem{definition}{Definition}[section]
\newtheorem{remark}{Remark}[section]

\numberwithin{equation}{section}
\numberwithin{figure}{section}
\numberwithin{table}{section}

\newcommand{\xx}{{\mathbf x}}

\begin{document}

\title{Normal form for GLT sequences, functions of normal GLT sequences, and spectral distribution of perturbed normal matrices}

\author{Giovanni Barbarino\\[-9pt]
\footnotesize Mathematics and Operational Research Unit, University of Mons, Belgium (giovanni.barbarino@umons.ac.be)\\
Carlo Garoni\thanks{Corresponding author}\\[-9pt]
\footnotesize Department of Mathematics, University of Rome Tor Vergata, Italy (garoni@mat.uniroma2.it)}
\date{}

\maketitle

\begin{abstract}
The theory of generalized locally Toeplitz (GLT) sequences is a powerful apparatus for computing the asymptotic spectral distribution of matrices $A_n$ arising from numerical discretizations of differential equations. Indeed, when the mesh fineness parameter $n$ tends to infinity, these matrices $A_n$ give rise to a sequence $\{A_n\}_n$, which often turns out to be a GLT sequence. In this paper, we extend the theory of GLT sequences in several directions: we show that every GLT sequence enjoys a normal form, we identify the spectral symbol of every GLT sequence formed by normal matrices, and we prove that, for every GLT sequence $\{A_n\}_n$ formed by normal matrices and every continuous function $f:\mathbb C\to\mathbb C$, the sequence $\{f(A_n)\}_n$ is again a GLT sequence whose spectral symbol is $f(\kappa)$, where $\kappa$ is the spectral symbol of $\{A_n\}_n$.
In addition, using the theory of GLT sequences, we prove a spectral distribution result for perturbed normal matrices.

\smallskip


\smallskip

\end{abstract}


\section{Introduction}

When a linear differential equation (DE) is discretized by a linear numerical method, the computation of the numerical solution reduces to solving a linear system $A_n\mathbf u_n=\mathbf f_n$, whose size $d_n$ increases with the mesh fineness parameter~$n$. What is often observed in practice is that $A_n$ enjoys an asymptotic spectral distribution in the limit of mesh refinement $n\to\infty$. More precisely, it often turns out that, for a large class of test functions~$F$,
\begin{equation*}
\lim_{n\to\infty}\frac1{d_n}\sum_{i=1}^{d_n}F(\lambda_i(A_n))=\frac1{\mu_k(D)}\int_DF(f(\xx)){\rm d}\xx,
\end{equation*}
where $\lambda_i(A_n)$, $i=1,\ldots,d_n$, are the eigenvalues of $A_n$, $\mu_k$ is the Lebesgue measure in $\mathbb R^k$, and $f:D\subset\mathbb R^k\to\mathbb C$. In this scenario, the function $f$ is referred to as the spectral symbol of the sequence $\{A_n\}_n$ and we write $\{A_n\}_n\sim_\lambda f$. We refer the reader to Remark~\ref{i.m.} for the informal meaning behind the spectral distribution $\{A_n\}_n\sim_\lambda f$ and to \cite[Chapter~1]{GLTbookI} for a list of practical uses of the spectral symbol $f$.

The theory of generalized locally Toeplitz (GLT) sequences is an apparatus---to the best of the authors' knowledge, the most powerful apparatus---for computing the spectral symbol $f$. Indeed, the sequence of discretization matrices $\{A_n\}_n$ turns out to be a GLT sequence for many DEs and numerical methods.
Nowadays, the main references for the theory of GLT sequences and related applications are the books \cite{GLTbookI,GLTbookII} and the papers \cite{barbarinoREDUCED,rGLT,GLTbookIII,GLTbookIV}. We therefore refer the reader to these works for a comprehensive treatment of the topic. 
For a more concise introduction to the subject, we recommend \cite{SbMath}.

In this paper, we extend the theory of GLT sequences in several directions. In particular:
\begin{itemize}[leftmargin=*,nolistsep]
	\item we show that every GLT sequence enjoys a normal form in the sense specified later on (see Theorem~\ref{normal});
	\item we identify the spectral symbol of every GLT sequence formed by normal matrices (see Theorem~\ref{GLT1-normal}); 
	\item we prove that, for every GLT sequence $\{A_n\}_n$ formed by normal matrices and every continuous function $f:\mathbb C\to\mathbb C$, the sequence $\{f(A_n)\}_n$ is again a GLT sequence whose spectral symbol is given by the composite function $f(\kappa)$, where $\kappa$ is the spectral symbol of $\{A_n\}_n$ (see Theorem~\ref{f(GLT)=GLT}).
\end{itemize}
In addition, using the theory of GLT sequences, we prove a spectral distribution result for perturbed normal matrices (see Theorem~\ref{perturbed-normals}).
In this regard, it is curious to note how a pure matrix analysis result apparently unrelated to GLT sequences can be proved using GLT sequences.
It is also worth pointing out that the proofs of the main results (Theorems~\ref{normal}--\ref{f(GLT)=GLT}) require some auxiliary results that are of interest also in themselves. We expressly refer to Theorems~\ref{diag2.5} and~\ref{N-cI}, which we have labeled as ``theorems'' in order to emphasize their importance over the other auxiliary results.

The paper is organized as follows. In Section~\ref{overview}, we overview the basics of the theory of GLT sequences. In Section~\ref{mr}, we state the main results. In Section~\ref{p}, we prove the main results.

\section{Overview of the theory of GLT sequences}\label{overview}
In this section, we overview the basics of the theory of GLT sequences. 

\subsection{Singular value and spectral distribution of a sequence of matrices}
Let $\mu_k$ be the Lebesgue measure in $\mathbb R^k$. Throughout this paper, all terminology from measure theory (such as ``measurable set'', ``measurable function'', ``a.e.'', etc.)\ always refers to the Lebesgue measure.
Let $C_c(\mathbb R)$ (resp., $C_c(\mathbb C)$) be the space of continuous complex-valued functions with bounded support defined on $\mathbb R$ (resp., $\mathbb C$). The singular values and eigenvalues of a matrix $A\in\mathbb C^{n\times n}$ are denoted by $\sigma_1(A),\ldots,\sigma_n(A)$ and $\lambda_1(A),\ldots,\lambda_n(A)$, respectively.

\begin{definition}[\textbf{singular value and spectral distribution of a sequence of matrices}]\label{dd}
Let $\{A_n\}_n$ be a sequence of matrices with $A_n$ of size $d_n\to\infty$, and let $f:D\subset\mathbb R^k\to\mathbb C$ be a measurable function defined on a set $D$ with $0<\mu_k(D)<\infty$.
\begin{itemize}[leftmargin=*,nolistsep]
	\item We say that $\{A_n\}_n$ has a spectral (or eigenvalue) distribution described by $f$, and we write $\{A_n\}_n\sim_\lambda f$, if
	\begin{equation}\label{sd}
	\lim_{n\to\infty}\frac1{d_n}\sum_{i=1}^{d_n}F(\lambda_i(A_n))=\frac1{\mu_k(D)}\int_DF(f(\xx)){\rm d}\xx,\qquad\forall\,F\in C_c(\mathbb C).
	\end{equation}
	In this case, $f$ is called the spectral (or eigenvalue) symbol of $\{A_n\}_n$.
	\item We say that $\{A_n\}_n$ has a singular value distribution described by $f$, and we write $\{A_n\}_n\sim_\sigma f$, if
	\begin{equation}\label{svd}
	\lim_{n\to\infty}\frac1{d_n}\sum_{i=1}^{d_n}F(\sigma_i(A_n))=\frac1{\mu_k(D)}\int_DF(|f(\xx)|){\rm d}\xx,\qquad\forall\,F\in C_c(\mathbb R).
	\end{equation}
	In this case, $f$ is called the singular value symbol of $\{A_n\}_n$.
\end{itemize}
\end{definition}


\begin{remark}\label{i.m.}
The informal meaning behind the spectral distribution \eqref{sd} is the following \cite[p.~46]{GLTbookI}: assuming that $f$ is continuous a.e., the eigenvalues of $A_n$, except possibly for $o(d_n)$ outliers, are approximately equal to the samples of $f$ over a uniform grid in the domain $D$ (for $n$ large enough).
A completely analogous meaning can be given for the singular value distribution \eqref{svd}.
We refer the reader to \cite[Theorems~2.1--2.2]{sdau}, \cite[Corollary~3.3]{DavideCalcolo2021}, \cite[Theorems~4.2 and~4.4]{Ftest}, \cite[Theorem~1.5]{maximum_norm} and \cite[Theorem~1.3]{Bottcher-ded-Grudsky} for a mathematical proof under certain assumptions of the informal meaning given here.
\end{remark}

\subsection{Special matrix-sequences}
In what follows, a matrix-sequence is a sequence of the form $\{A_n\}_n$, where $A_n\in\mathbb C^{n\times n}$ for every~$n$.
If $A\in\mathbb C^{n\times n}$ and $1\le p\le\infty$, we denote by $\|A\|_p$ the Schatten $p$-norm of $A$, i.e., the $p$-norm of the vector $(\sigma_1(A),\ldots,\sigma_n(A))$ formed by the singular values of $A$. The Schatten $2$-norm $\|A\|_2$ coincides with the Frobenius norm. The Schatten $\infty$-norm $\|A\|_\infty$ coincides with the classical spectral (or Euclidean) norm and is preferably denoted by $\|A\|$. For more on Schatten $p$-norms, see \cite{Bhatia}.

\subsubsection{Zero-distributed sequences}
A matrix-sequence $\{Z_n\}_n$ such that $Z_n=R_n+N_n$ with
\[ \lim_{n\to\infty}n^{-1}\mathrm{rank}(R_n)=\lim_{n\to\infty}\|N_n\|=0 \]
is referred to as a zero-distributed sequence. It can be shown that a matrix-sequence $\{Z_n\}_n$ is zero-distributed if and only if $\{Z_n\}_n\sim_\sigma0$, i.e., 
\[ \lim_{n\to\infty}\frac1{n}\sum_{i=1}^nF(\sigma_i(Z_n))=F(0),\qquad\forall\,F\in C_c(\mathbb R); \]
see \cite[Chapter~3]{GLTbookI}. 

\subsubsection{Diagonal sampling sequences}
If $n\in\mathbb N$ and $a:[0,1]\to\mathbb C$, the $n$th diagonal sampling matrix generated by $a$ is the $n\times n$ diagonal matrix given by
\begin{equation*}
D_n(a)=\mathop{\rm diag}_{i=1,\ldots,n}a\Bigl(\frac in\Bigr).
\end{equation*}
$\{D_n(a)\}_n$ is called the diagonal sampling sequence generated by $a$.

\subsubsection{Toeplitz sequences}
If $n\in\mathbb N$ and $f:[-\pi,\pi]\to\mathbb C$ is a function in $L^1([-\pi,\pi])$, the $n$th Toeplitz matrix generated by $f$ is the $n\times n$ matrix
\[ T_n(f)=[f_{i-j}]_{i,j=1}^n, \] 
where the numbers $f_k$ are the Fourier coefficients of $f$:
\begin{equation}\label{fF}
f_k=\frac1{2\pi}\int_{-\pi}^\pi f(\theta){\rm e}^{-{\rm i}k\theta}{\rm d}\theta,\qquad k\in\mathbb Z.
\end{equation}
$\{T_n(f)\}_n$ is called the Toeplitz sequence generated by $f$. For more details on Toeplitz matrices generated by a function $f$, see \cite{BoSi} and \cite[Chapter~6]{GLTbookI}.

\subsection{Approximating classes of sequences}

\begin{definition}[\textbf{approximating class of sequences}]
Let $\{A_n\}_n$ be a matrix-sequence and let $\{\{B_{n,m}\}_n\}_m$ be a sequence of matrix-sequences. We say that $\{\{B_{n,m}\}_n\}_m$ is an approximating class of sequences (a.c.s.)\ for $\{A_n\}_n$, and we write $\{B_{n,m}\}_n\xrightarrow{\rm a.c.s.}\{A_n\}_n$, if the following condition is met: for every $m$ there exists $n_m$ such that, for $n\ge n_m$,
\begin{equation*}
A_n=B_{n,m}+R_{n,m}+N_{n,m},\qquad{\rm rank}(R_{n,m})\le c(m)n,\qquad\|N_{n,m}\|\le\omega(m),
\end{equation*}
where $n_m,\,c(m),\,\omega(m)$ depend only on $m$ and $\lim_{m\to\infty}c(m)=\lim_{m\to\infty}\omega(m)=0$.
\end{definition}

Roughly speaking, $\{\{B_{n,m}\}_n\}_m$ is an a.c.s.\ for $\{A_n\}_n$ if, for all sufficiently large $m$, the sequence $\{B_{n,m}\}_n$ approximates $\{A_n\}_n$ in the sense that $A_n$ is eventually equal to $B_{n,m}$ plus a small-rank matrix (with respect to the matrix size $n$) plus a small-norm matrix. 
It turns out that the notion of a.c.s.\ is a notion of convergence in the space of matrix-sequences
\[ \mathscr E=\{\{A_n\}_n:\{A_n\}_n\mbox{ is a matrix-sequence}\}. \]
More precisely, for every $A\in\mathbb C^{n\times n}$ let 
\begin{align}
p(A)&=\inf\biggl\{\frac{{\rm rank}(R)}{n}+\|N\|:\quad R,N\in\mathbb C^{n\times n},\quad R+N=A\biggr\}=\min_{i=0,\ldots,n}\Bigl(\frac in+\sigma_{i+1}(A)\Bigr),\label{p(A)-min}
\end{align}
where $\sigma_1(A)\ge\ldots\ge\sigma_n(A)$ and $\sigma_{n+1}(A)=0$ by convention. For the second equation in \eqref{p(A)-min}, see \cite[Section~5.2.2]{GLTbookI}.
Set
\begin{align*}
p_{\rm a.c.s.}(\{A_n\}_n)&=\limsup_{n\to\infty}\,p(A_n),\qquad \{A_n\}_n\in\mathscr E,\\
d_{\rm a.c.s.}(\{A_n\}_n,\{B_n\}_n)&=p_{\rm a.c.s.}(\{A_n-B_n\}_n),\qquad\{A_n\}_n,\{B_n\}_n\in\mathscr E.
\end{align*}
It is proved in \cite[Chapter~5]{GLTbookI} that $d_{\rm a.c.s.}$ is a distance on $\mathscr E$ which turns $\mathscr E$ into a pseudometric space $(\mathscr E,d_{\rm a.c.s.})$ where the statement ``$\{\{B_{n,m}\}_n\}_m$ converges to $\{A_n\}_n$'' is equivalent to ``$\{\{B_{n,m}\}_n\}_m$ is an a.c.s.\ for $\{A_n\}_n$''. In other words, we have
\begin{equation*}
\{B_{n,m}\}_n\xrightarrow{\rm a.c.s.}\{A_n\}_n\quad\iff\quad\lim_{m\to\infty}d_{\rm a.c.s.}(\{B_{n,m}\}_n,\{A_n\}_n)=0.
\end{equation*}
We remark that the a.c.s.\ distance $d_{\rm a.c.s.}$ is a pseudometric but not a metric, because $d_{\rm a.c.s.}(\{A_n\}_n,\{B_n\}_n)$ can be zero even if $\{A_n\}_n$ is not equal to $\{B_n\}_n$. In fact, we have 
\begin{equation}\label{acs-equiv}
d_{\rm a.c.s.}(\{A_n\}_n,\{B_n\}_n)=0\quad\iff\quad\{A_n-B_n\}_n\mbox{ is zero-distributed};
\end{equation}
see \cite[Section~5.2.1]{GLTbookI}. Whenever one of the equivalent conditions in \eqref{acs-equiv} is satisfied, we say that $\{A_n\}_n$ is a.c.s.-equivalent to $\{B_n\}_n$ and we write $\{A_n\}_n\equiv_{\rm a.c.s.}\{B_n\}_n$. 
We refer the reader to \cite[Chapter~5]{GLTbookI} and \cite{d-complete,ela} for deeper insights into the a.c.s.\ convergence.

\subsection{GLT sequences}

\begin{definition}[\textbf{GLT sequence}]
Let $\{A_n\}_n$ be a matrix-sequence and let $\kappa:[0,1]\times[-\pi,\pi]\to\mathbb C$ be measurable. We say that $\{A_n\}_n$ is a GLT sequence with symbol $\kappa$, and we write $\{A_n\}_n\sim_{\rm GLT}\kappa$, if there exist functions $a_{i,m}$, $f_{i,m}$, $i=1,\ldots,N_m$, with $a_{i,m}:[0,1]\to\mathbb C$ continuous a.e.\ and $f_{i,m}:[-\pi,\pi]\to\mathbb C$ in $L^1([-\pi,\pi])$, such that
\begin{itemize}[nolistsep,leftmargin=*]
	\item $\kappa_m(x,\theta)=\sum_{i=1}^{N_m}a_{i,m}(x)f_{i,m}(\theta)\to\kappa(x,\theta)$ in measure over $[0,1]\times[-\pi,\pi]$, \vspace{2pt}
	\item $\{A_{n,m}\}_n=\bigl\{\sum_{i=1}^{N_m}D_n(a_{i,m})T_n(f_{i,m})\bigr\}_n\xrightarrow{\rm a.c.s.}\{A_n\}_n$.
\end{itemize}
\end{definition}
We report below the properties of GLT sequences that we need in this paper. The corresponding proofs can be found in \cite[Chapter~8]{GLTbookI}. 
\begin{enumerate}[leftmargin=39pt,nolistsep]
	\item[\textbf{GLT\,1.}] If $\{A_n\}_n\sim_{\rm GLT}\kappa$ then $\{A_n\}_n\sim_\sigma\kappa$. If $\{A_n\}_n\sim_{\rm GLT}\kappa$ and the matrices $A_n$ are Hermitian then $\{A_n\}_n\sim_\lambda\kappa$.
	\item[\textbf{GLT\,2.}] We have
	\begin{itemize}[nolistsep,leftmargin=*]
		\item $\{T_n(f)\}_n\sim_{\rm GLT}\kappa(x,\theta)=f(\theta)$ if $f\in L^1([-\pi,\pi])$,
		\item $\{D_n(a)\}_n\sim_{\rm GLT}\kappa(x,\theta)=a(x)$ if $a:[0,1]\to\mathbb C$ is continuous a.e.,
		\item $\{Z_n\}_n\sim_{\rm GLT}\kappa(x,\theta)=0$ if and only if $\{Z_n\}_n$ is zero-distributed.
	\end{itemize}
	\item[\textbf{GLT\,3.}] If $\{A_n\}_n\sim_{\rm GLT}\kappa$ and $\{B_n\}_n\sim_{\rm GLT}\xi$ then
	\begin{itemize}[nolistsep,leftmargin=*]
		\item $\{A_n^*\}_n\sim_{\rm GLT}\overline\kappa$,
		\item $\{\alpha A_n+\beta B_n\}_n\sim_{\rm GLT}\alpha\kappa+\beta\xi$ for all $\alpha,\beta\in\mathbb C$,
		\item $\{A_nB_n\}_n\sim_{\rm GLT}\kappa\xi$,
	\end{itemize}
	\item[\textbf{GLT\,4.}] If $\{A_n\}_n\sim_{\rm GLT}\kappa$ and the matrices $A_n$ are Hermitian then $\{f(A_n)\}_n\sim_{\rm GLT}f(\kappa)$ for every continuous function $f:\mathbb C\to\mathbb C$.
	\item[\textbf{GLT\,5.}] $\{A_n\}_n\sim_{\rm GLT}\kappa$ if and only if there exist GLT sequences $\{B_{n,m}\}_n\sim_{\rm GLT}\kappa_m$ such that $\{B_{n,m}\}_n\xrightarrow{\rm a.c.s.}\{A_n\}_n$ and $\kappa_m\to\kappa$ in measure.
	\item[\textbf{GLT\,6.}] Suppose $\{A_n\}_n\sim_{\rm GLT}\kappa$ and $\{B_{n,m}\}_n\sim_{\rm GLT}\kappa_m$. Then, $\{B_{n,m}\}_n\xrightarrow{\rm a.c.s.}\{A_n\}_n$ if and only if $\kappa_m\to\kappa$ in measure.
\end{enumerate}

\section{Main results}\label{mr}

We have now collected all the necessary ingredients to state our main results. Throughout this paper, we denote by $O_n$ and $I_n$ the $n\times n$ zero matrix and the $n\times n$ identity matrix, respectively. Moreover, we denote by
\begin{equation}\label{Fn.expr}
F_n=\frac1{\sqrt n}\bigl[{\rm e}^{-2\pi{\rm i}jk/n}\bigr]_{j,k=0}^{n-1}=\frac1{\sqrt n}\bigl[{\rm e}^{-2\pi{\rm i}(j-1)(k-1)/n}\bigr]_{j,k=1}^n
\end{equation}
the unitary discrete Fourier transform of size $n$. For every $m,n\in\mathbb N$ with $m\le n$, we define the normal matrix
\begin{equation}\label{Qnm}
Q_{n,m}=\left[\:
\begin{array}{ccccc}
\cline{1-1}
\multicolumn{1}{|c|}{F_{\lfloor n/m\rfloor}} &&&& \\
\cline{1-2} 
& \multicolumn{1}{|c|}{F_{\lfloor n/m\rfloor}} &&&\\
\cline{2-2}
&& \ddots &&\\
\cline{4-4}
&&&  \multicolumn{1}{|c|}{F_{\lfloor n/m\rfloor}}
&\\
\cline{4-5}
&&&& \multicolumn{1}{|c|}{I_{n\,{\rm mod}\,m}}
\\\cline{5-5}
\end{array}
\:\right],
\end{equation}
where it is understood that the last identity block $I_{n\,{\rm mod}\,m}$ is not present if $n$ is a multiple on $m$. 
A matrix-sequence formed by diagonal (resp., unitary, normal) matrices is referred to as a diagonal (resp., unitary, normal) sequence.
The first main result of this paper (Theorem~\ref{normal}) delivers a normal form for every GLT sequence.

\begin{theorem}[\textbf{normal form for GLT sequences}]\label{normal}
There exists a unitary sequence $\{Q_n\}_n$ with the following property.
\begin{quote}
For every GLT sequence $\{A_n\}_n\sim_{\rm GLT}\kappa$ there exists a diagonal sequence $\{D_n\}_n$ such that $\{A_n\}_n\equiv_{\rm a.c.s.}\{Q_nD_nQ_n^*\}_n$, $\{D_n\}_n\sim_\lambda\kappa$, and $D_n$ is real for every $n$ if $\kappa$ is real. 
\end{quote}
Moreover, if $\{m(n)\}_n$ is any sequence such that $m(n)\le n$, $m(n)\to\infty$ and $m(n)=o(n)$ as $n\to\infty$, then $Q_n$ can be chosen as $Q_n=Q_{n,m(n)}$.
\end{theorem}

As we shall see, the matrices $D_n$ are built in order to contain an ordered approximation of the samples of the symbol $\kappa$ over $[0,1]\times[-\pi,\pi]$.
Note that a unitary transformation does not change the eigenvalues, so $\{Q_nD_nQ_n^*\}_n\sim_\lambda\kappa$.
We refer to the sequence $\{Q_nD_nQ_n^*\}_n$ as a normal form for the GLT sequence $\{A_n\}_n$. 

The second main result of this paper (Theorem~\ref{perturbed-normals}) provides the spectral distribution of matrix-sequences of the form $\{X_n+Y_n\}_n$, where $\{X_n\}_n$ is a normal sequence such that $\{X_n\}_n\sim_\lambda f$ and $Y_n$ is a ``small'' perturbation of $X_n$. Throughout this paper, we use the natural convention $1/\infty=0$.

\begin{theorem}[\textbf{spectral distribution of perturbed normal sequences}]\label{perturbed-normals}
Let $\{X_n\}_n$ be a normal sequence such that $\{X_n\}_n\sim_\lambda f$ and let $\{Y_n\}_n$ be a matrix-sequence. Suppose that one of the following conditions is met.
\begin{enumerate}[nolistsep,leftmargin=*]
	\item $\|Y_n\|_p=o(1)$ for some $1\le p\le 2$.
	\item $\|Y_n\|_p=o(n^{2/p-1})$ for some $2\le p\le\infty$.
	\item $\{Y_n\}_n$ is zero-distributed and the matrices $X_n+Y_n$ are normal.
\end{enumerate}
Then $\{X_n+Y_n\}_n\sim_\lambda f$.
\end{theorem}


The last two main results of this paper (Theorems~\ref{GLT1-normal} and~\ref{f(GLT)=GLT}) generalize properties {\bf GLT\,1} and {\bf GLT\,4} to the case where the matrices $A_n$ are normal.

\begin{theorem}[\textbf{spectral distribution of normal GLT sequences}]\label{GLT1-normal}
If $\{A_n\}_n\sim_{\rm GLT}\kappa$ and the matrices $A_n$ are normal then $\{A_n\}_n\sim_\lambda\kappa$.
\end{theorem}

\begin{theorem}[\textbf{functions of normal GLT sequences}]\label{f(GLT)=GLT}
If $\{A_n\}_n\sim_{\rm GLT}\kappa$ and the matrices $A_n$ are normal then $\{f(A_n)\}_n\sim_{\rm GLT}f(\kappa)$ for every continuous function $f:\mathbb C\to\mathbb C$.
\end{theorem}

Note that $f(A_n)$ is normal whenever $A_n$ is normal. Thus, under the hypotheses of Theorem~\ref{f(GLT)=GLT}, we have $\{f(A_n)\}_n\sim_\lambda f(\kappa)$ for every continuous function $f:\mathbb C\to\mathbb C$ by Theorem~\ref{GLT1-normal}.

\section{Proofs of the main results}\label{p}

\subsection{Proof of Theorem~\ref{normal}}\label{p1}
The proof of Theorem~\ref{normal} requires several preliminary results.
In particular, Theorem~\ref{diag2.5} is analogous to the main a.c.s.\ approximation results \cite[Corollaries~5.1 and~5.2]{GLTbookI} and it is therefore of interest also in itself.

\subsubsection{Diagonal sequences}
If $A=U\Sigma V$ is an SVD of $A\in\mathbb C^{n\times n}$ (recall the the singular values $\sigma_1(A)\ge\ldots\ge\sigma_n(A)$ are arranged in decreasing order in $\Sigma$ by definition of SVD), then we set, for every $i=0,\ldots,n$, 
\begin{align*}
\hat A^{(i)}&=U\hat\Sigma^{(i)}V=U\,{\rm diag}(\sigma_1(A),\ldots,\sigma_i(A),0,\ldots,0)\,V,\\
\tilde A^{(i)}&=U\tilde\Sigma^{(i)}V=U\,{\rm diag}(0,\ldots,0,\sigma_{i+1}(A),\ldots,\sigma_n(A))\,V.
\end{align*}
For every $A\in\mathbb C^{n\times n}$ and every SVD $A=U\Sigma V$, the scalar quantity $p(A)$ in \eqref{p(A)-min} satisfies
\[ p(A)=\min_{i=0,\ldots,n}\Bigl(\frac in+\sigma_{i+1}(A)\Bigr)=\min_{i=0,\ldots,n}\Bigl(\frac{{\rm rank}(\hat A^{(i)})}n+\|\tilde A^{(i)}\|\Bigr), \]
because the minimum in both sides of the last equality is reached for $i\in\{0,\ldots,{\rm rank}(A)\}$, since for $i\ge{\rm rank}(A)$ we have $\|\tilde A^{(i)}\|=\sigma_{i+1}(A)=0$.

\begin{lemma}\label{diag2.2}
Let $\{A_n\}_n$ and $\{B_{n,m}\}_n$ be matrix-sequences, let $A_n-B_{n,m}=U_{n,m}\Sigma_{n,m}V_{n,m}$ be an SVD of $A_n-B_{n,m}$, define
\begin{alignat*}{3}
R_{n,m}^{(i)}&=\widehat{(A_n-B_{n,m})}^{(i)}, &\quad i&=0,\ldots,n,\\
N_{n,m}^{(i)}&=\widetilde{(A_n-B_{n,m})}^{(i)}, &\quad i&=0,\ldots,n,
\end{alignat*}
and let $i_{n,m}\in\{0,\ldots,n\}$ be the index such that
\[ p(A_n-B_{n,m})=\min_{i=0,\ldots,n}\Bigl(\frac{{\rm rank}(R_{n,m}^{(i)})}n+\|N_{n,m}^{(i)}\|\Bigr)=\frac{{\rm rank}(R_{n,m}^{(i_{n,m})})}{n}+\|N_{n,m}^{(i_{n,m})}\|.  \]
Then, $\{B_{n,m}\}_n\xrightarrow{\rm a.c.s.}\{A_n\}_n$ if and only if for every $m$ there exists $n_m$ such that, for $n\ge n_m$,
\[ {\rm rank}(R_{n,m}^{(i_{n,m})})\le c(m)n,\qquad\|N_{n,m}^{(i_{n,m})}\|\le\omega(m), \]
where $n_m,c(m),\omega(m)$ depend only on $m$ and $\lim_{m\to\infty}c(m)=\lim_{m\to\infty}\omega(m)=0$.
\end{lemma}
\begin{proof}
We have
\begin{align*}
\{B_{n,m}\}_n\xrightarrow{\rm a.c.s.}\{A_n\}_n&\quad\iff\quad \lim_{m\to\infty}d_{\rm a.c.s.}(\{A_n\}_n,\{B_{n,m}\}_n)=0\\
&\quad\iff\quad\lim_{m\to\infty}\limsup_{n\to\infty}p(A_n-B_{n,m})=0\\
&\quad\iff\quad\lim_{m\to\infty}\limsup_{n\to\infty}\Bigl(\frac{{\rm rank}(R_{n,m}^{(i_{n,m})})}{n}+\|N_{n,m}^{(i_{n,m})}\|\Bigr)=0\\
&\quad\iff\quad\lim_{m\to\infty}\limsup_{n\to\infty}\frac{{\rm rank}(R_{n,m}^{(i_{n,m})})}{n}=0\quad\wedge\quad\lim_{m\to\infty}\limsup_{n\to\infty}\|N_{n,m}^{(i_{n,m})}\|=0.
\end{align*}
The last two conditions are satisfied if and only if for every $m$ there exists $n_m$ such that, for $n\ge n_m$,
\[ {\rm rank}(R_{n,m}^{(i_{n,m})})\le c(m)n,\qquad\|N_{n,m}^{(i_{n,m})}\|\le\omega(m), \]
where $n_m,c(m),\omega(m)$ depend only on $m$ and $\lim_{m\to\infty}c(m)=\lim_{m\to\infty}\omega(m)=0$.
\end{proof}

\begin{corollary}\label{diag2.2'}
Let $\{A_n\}_n$ and $\{B_{n,m}\}_n$ be diagonal sequences, and suppose that $\{B_{n,m}\}_n\xrightarrow{\rm a.c.s.}\{A_n\}_n$. Then, for every $m$ there exists $n_m$ such that, for $n\ge n_m$,
\[ A_n=B_{n,m}+R_{n,m}+N_{n,m},\qquad{\rm rank}(R_{n,m})\le c(m)n,\qquad\|N_{n,m}\|\le\omega(m), \]
where $R_{n,m},N_{n,m}$ are diagonal matrices, $n_m,c(m),\omega(m)$ depend only on $m$, and $\lim_{m\to\infty}c(m)\hspace{-0.5pt}=\hspace{-0.5pt}\lim_{m\to\infty}\omega(m)\hspace{-0.5pt}=\hspace{-0.5pt}0$.
\end{corollary}
\begin{proof}
Since $A_n-B_{n,m}$ is diagonal, there exists an SVD $A_n-B_{n,m}=U_{n,m}\Sigma_{n,m}V_{n,m}$ with $U_{n,m}$ and $V_{n,m}$ diagonal.
By applying Lemma~\ref{diag2.2} with this SVD, we conclude that for every $m$ there exists $n_m$ such that, for $n\ge n_m$,
\[ A_n=B_{n,m}+R_{n,m}^{(i_{n,m})}+N_{n,m}^{(i_{n,m})},\qquad{\rm rank}(R_{n,m}^{(i_{n,m})})\le c(m)n,\qquad\|N_{n,m}^{(i_{n,m})}\|\le\omega(m), \]
where $n_m,c(m),\omega(m)$ depend only on $m$ and $\lim_{m\to\infty}c(m)=\lim_{m\to\infty}\omega(m)=0$. Since $R_{n,m}^{(i_{n,m})}$ and $N_{n,m}^{(i_{n,m})}$ are diagonal by definition (and by the fact that $U_{n,m}$ and $V_{n,m}$ are diagonal), the result is proved.
\end{proof}


\begin{theorem}\label{diag2.5}
Let $\{A_n\}_n$, $\{B_{n,m}\}_n$ be matrix-sequences formed by diagonal matrices and let $f,f_m:D\subset\mathbb R^k\to\mathbb C$ be measurable functions defined on a set $D$ with $0<\mu_k(D)<\infty$. Suppose that
\begin{enumerate}[nolistsep,leftmargin=*]
	\item $\{B_{n,m}\}_n\sim_\lambda f_m$ for every $m$,
	\item $\{B_{n,m}\}_n\xrightarrow{\rm a.c.s.}\{A_n\}_n$,
	\item $f_m\to f$ in measure.
\end{enumerate}
Then $\{A_n\}_n\sim_\lambda f$.
\end{theorem}
\begin{proof}
Let $F\in C_c(\mathbb C)$. For all $n,m$ we have
\begin{align}
\left|\frac1n\sum_{j=1}^nF(\lambda_j(A_n))-\frac1{\mu_k(D)}\int_DF(f(\xx)){\rm d}\xx\right|&\le\left|\frac1n\sum_{j=1}^nF(\lambda_j(A_n))-\frac1n\sum_{j=1}^nF(\lambda_j(B_{n,m}))\right|\notag\\
&\qquad+\left|\frac1n\sum_{j=1}^nF(\lambda_j(B_{n,m}))-\frac1{\mu_k(D)}\int_DF(f_m(\xx)){\rm d}\xx\right|\notag\\
&\qquad+\left|\frac1{\mu_k(D)}\int_DF(f_m(\xx)){\rm d}\xx-\frac1{\mu_k(D)}\int_DF(f(\xx)){\rm d}\xx\right|.\label{chain<=}
\end{align}
The second term in the right-hand side of \eqref{chain<=} tends to 0 as $n\to\infty$ because $\{B_{n,m}\}_n\sim_\lambda f_m$ for every $m$ by hypothesis. The third term in the right-hand side of \eqref{chain<=} tends to 0 as $m\to\infty$ because $f_m\to f$ in measure by hypothesis, which implies that $\int_DF(f_m(\xx)){\rm d}\xx\to\int_DF(f(\xx)){\rm d}\xx$; see, e.g., \cite[Lemma~2.5]{GLTbookI}. We show that
\begin{equation}\label{limlimsup}
\lim_{m\to\infty}\limsup_{n\to\infty}\left|\frac1n\sum_{j=1}^nF(\lambda_j(A_n))-\frac1n\sum_{j=1}^nF(\lambda_j(B_{n,m}))\right|=0,
\end{equation}
after which the thesis follows by passing first to the $\limsup_{n\to\infty}$ and then to the $\lim_{m\to\infty}$ in \eqref{chain<=}. To prove \eqref{limlimsup}, recall that $\{B_{n,m}\}_n\xrightarrow{\rm a.c.s.}\{A_n\}_n$ and $A_n,B_{n,m}$ are diagonal. Hence, by Corollary~\ref{diag2.2'}, for every $m$ there exists $n_m$ such that, for $n\ge n_m$,
\[ A_n=B_{n,m}+R_{n,m}+N_{n,m},\qquad{\rm rank}(R_{n,m})\le c(m)n,\qquad\|N_{n,m}\|\le\omega(m), \]
where $R_{n,m},N_{n,m}$ are diagonal and $\lim_{m\to\infty}c(m)=\lim_{m\to\infty}\omega(m)=0$. It follows that
\[ \#\Bigl\{j\in\{1,\ldots,n\}:|\lambda_j(A_n)-\lambda_j(B_{n,m})|>\|N_{n,m}\|\Bigr\}\le{\rm rank}(R_{n,m}). \]
Indeed, since
\[ \lambda_j(A_n)-\lambda_j(B_{n,m})=\lambda_j(A_n-B_{n,m})=\lambda_j(R_{n,m}+N_{n,m})=\lambda_j(R_{n,m})+\lambda_j(N_{n,m}), \]
if we have $|\lambda_j(A_n)-\lambda_j(B_{n,m})|>\|N_{n,m}\|$ then we also have $\lambda_j(R_{n,m})\ne0$, because $|\lambda_j(N_{n,m})|\le\|N_{n,m}\|$. 
Thus, for every $m$ and every $n\ge n_m$, if we denote by $\omega_F$ the modulus of continuity of $F$, we have
\begin{align*}
&\left|\frac1n\sum_{j=1}^nF(\lambda_j(A_n))-\frac1n\sum_{j=1}^nF(\lambda_j(B_{n,m}))\right|\\
&\le\frac1n\sum_{\substack{j\in\{1,\ldots,n\}:\\|\lambda_j(A_n)-\lambda_j(B_{n,m})|\le\|N_{n,m}\|}}|F(\lambda_j(A_n))-F(\lambda_j(B_{n,m}))|+\frac1n\sum_{\substack{j\in\{1,\ldots,n\}:\\|\lambda_j(A_n)-\lambda_j(B_{n,m})|>\|N_{n,m}\|}}|F(\lambda_j(A_n))-F(\lambda_j(B_{n,m}))|\\
&\le\omega_F(\|N_{n,m}\|)+2\|F\|_\infty\frac{{\rm rank}(R_{n,m})}{n}\\
&\le\omega_F(\omega(m))+2\|F\|_\infty c(m),
\end{align*}
which tends to $0$ as $m\to\infty$. 
\end{proof}

\subsubsection{Circulant sequences}
If $n\in\mathbb N$ and $f:[-\pi,\pi]\to\mathbb C$ is a function in $L^1([-\pi,\pi])$, the $n$th circulant matrix generated by $f$ is the $n\times n$ matrix
\[ C_n(f)=\sum_{k=-(n-1)}^{n-1}f_kC_n^k, \]
where the numbers $f_k$ are the Fourier coefficient of $f$ defined in \eqref{fF} and
\begin{equation*}
C_n=\left[\begin{array}{ccccc}
0 & \ \ & \ \ & \ \ 1\\
1 & \ \ \ddots & \ \ & \ \ \\
& \ \ \ddots & \ \ \ddots & \ \ \\
& \ \ & \ \ 1 & \ \ 0 
\end{array}\right].
\end{equation*}
$C_n$ is called the generator of circulant matrices of size $n$ and $\{C_n(f)\}_n$ is called the circulant sequence generated by $f$. For the proof of the next theorem, see \cite[Section~6.4]{GLTbookI}. For more details on circulant matrices, see \cite{Davis}.
In what follows, if $n\in\mathbb N$ and $g:[0,2\pi]\to\mathbb C$, we denote by $\Delta_n(g)$ the diagonal sampling matrix
\[ \Delta_n(g)=\mathop{\rm diag}_{j=0,\ldots,\,n-1}g\Bigl(\frac{2\pi j}n\Bigr). \]
Moreover, if $t\in\mathbb N$ and $f\in L^1([-\pi,\pi])$, we denote by $f^{(t)}$ the truncated Fourier series
\[ f^{(t)}(\theta)=\sum_{k=-t}^{t}f_k{\rm e}^{{\rm i}k\theta}. \]

\begin{theorem}\label{sigmaC}
Let $n\in\mathbb N$ and $f\in L^1([-\pi,\pi])$. Then,
\begin{equation*}
C_n(f)=F_n\Delta_n(f^{(n-1)})F_n^*,
\end{equation*}
where $F_n$ is defined in \eqref{Fn.expr}. In particular, $C_n(f)$ is a normal matrix whose spectrum is given by
\[ \biggl\{f^{(n-1)}\Bigl(\frac{2\pi j}n\Bigr):\,j=0,\ldots,n-1\biggr\}. \]
\end{theorem}

\subsubsection{LT operator}

If $X,Y$ are matrices of any size, say $X\in\mathbb C^{m_1\times m_2}$ and $Y\in\mathbb C^{\ell_1\times \ell_2}$, the tensor (Kronecker) product of $X$ and $Y$ is the $m_1\ell_1\times m_2\ell_2$ matrix defined by
\[ X\otimes Y=\bigl[x_{ij}Y\bigr]_{\substack{i=1,\ldots,m_1\\ j=1,\ldots,m_2}}=\begin{bmatrix}
x_{11}Y & \cdots & x_{1m_2}Y\\
\vdots & & \vdots\\
x_{m_11}Y & \cdots & x_{m_1m_2}Y
\end{bmatrix}, \]
and the direct sum of $X$ and $Y$ is the $(m_1+\ell_1)\times(m_2+\ell_2)$ matrix defined by
\[ X\oplus Y=\mathop{\rm diag}(X,Y)=\begin{bmatrix}X & O\\O & Y\end{bmatrix}. \]

\begin{definition}[\textbf{locally Toeplitz operator}]\label{LTop}
Let $m,n\in\mathbb N$ with $m\le n$, let $a:[0,1]\to\mathbb C$ and let $f\in L^1([-\pi,\pi])$. The locally Toeplitz (LT) operator is defined as the following $n\times n$ matrix:
\begin{align*}
LT_n^m(a,f)&=D_m(a)\otimes T_{\left\lfloor n/m\right\rfloor}(f)\:\oplus\:O_{n\,{\rm mod}\,m}=\mathop{\rm diag}_{i=1,\ldots,m}\biggl[a\Bigl(\frac im\Bigr)T_{\left\lfloor n/m\right\rfloor}(f)\biggr]\:\oplus\:O_{n\,{\rm mod}\,m}\\
&=\left[\:
\begin{array}{ccccc}
\cline{1-1}
\multicolumn{1}{|c|}{a(\frac1m)T_{\lfloor n/m\rfloor}(f)} &&&& \\
\cline{1-2} 
& \multicolumn{1}{|c|}{a(\frac2m)T_{\lfloor n/m\rfloor}(f)} &&&\\
\cline{2-2}
&& \ddots &&\\
\cline{4-4}
&&&  \multicolumn{1}{|c|}{a(1)T_{\lfloor n/m\rfloor}(f)}
&\\
\cline{4-5}
&&&& \multicolumn{1}{|c|}{O_{n\,{\rm mod}\,m}}
\\\cline{5-5}
\end{array}
\:\right].
\end{align*}
It is understood that the zero block $O_{n\,{\rm mod}\,m}$ is not present when $n$ is a multiple of $m$. Moreover, here and in what follows, the tensor product operation $\otimes$ is always applied before the direct sum $\oplus$, exactly as in the case of numbers, where multiplication is always applied before addition.
\end{definition}

For our purposes, we are interested only in $LT_n^m(a,f)$ with $m=m(n)$ such that $m(n)\to\infty$ and $m(n)=o(n)$ as $n\to\infty$.

\begin{lemma}\label{LTacsDT}
Let $a:[0,1]\to\mathbb C$ be continuous and let $f$ be a trigonometric polynomial. Then,
\[ \{LT_n^{m(n)}(a,f)\}_n\equiv_{\rm a.c.s.}\{D_n(a)T_n(f)\}_n \]
for every $m(n)\to\infty$ such that $m(n)=o(n)$.
\end{lemma}
\begin{proof}
Throughout this proof, $m=m(n)$ and $m(n)$ has the properties specified in the statement.
We have to show that $\{Z_n\}_n$ is zero-distributed, where
\[ Z_n = LT_n^m(a,f) - D_n(a)T_n(f). \]
The matrix $Z_n$ is banded with bandwidth $2r+1$, where $r$ is the degree of $f$. Let $R_n$ be the matrix whose non-zero entries are the non-zero entries of $Z_n$ lying outside the diagonal blocks of $LT_n^m(a,f)$ or in its last zero block of size $n\,{\rm mod}\,m$. 
The matrix $R_n$ has at most $2rm + n\,{\rm mod}\,m=o(n)$ non-zero rows. 
Hence, ${\rm rank}(R_n) = o(n)$.
Now, let $N_n=Z_n-R_n$. The only non-zero entries of $N_n$ lie inside the diagonal blocks of $LT_n^m(a,f)$ (excluding the last zero block). Each non-zero entry of $N_n$ is of the form $(a(x)-a(y))\tau$, where $x,y\in[0,1]$ and $\tau$ is a coefficient of $T_n(f)$. Moreover, the following properties hold.
\begin{itemize}[nolistsep,leftmargin=*]
	\item The distance between $x$ and $y$ is small. Indeed, if we are in the $i$th block with $1\le i\le m$, then 
	\[ x = \frac im,\qquad ny\in \Bigl[(i-1)\Big\lfloor \frac nm\Big\rfloor + 1 ,\,i\Big\lfloor \frac nm\Big\rfloor\Bigr]. \] 
	Since
	\[ \frac im \ge \frac in\Bigl\lfloor\frac nm\Bigr\rfloor, \]
	we have
	\[ |x-y| \le \frac im  - \frac{i-1}n\Bigl\lfloor \frac nm\Bigr\rfloor \le \frac im-\frac{i-1}n\Bigl(\frac nm-1\Bigr)=\frac1m+\frac{i-1}n\le\frac1m+\frac mn. \]
	It follows that $|a(x) - a(y)| \le \omega_a(|x-y|)\le\omega_a(1/m+m/n)$.
	\item $|\tau|\le M$, where $M=\max_{k=-r,\ldots,r}|f_k|$.
\end{itemize}
Since $N_n$ is banded with bandwidth $2r+1$, we conclude that
\[ \|N_n\|\le\sqrt{\|N_n\|_1\|N_n\|_\infty}\le(2r+1)\,\omega_a\Bigl(\frac1m+\frac mn\Bigr)\,M=o(1). \]
In conclusion, $Z_n=R_n+N_n$ with ${\rm rank}(R_n)=o(n)$ and $\|N_n\|=o(1)$, which means that $\{Z_n\}_n$ is zero-distributed.
\end{proof}

\subsubsection{LC operator}\label{lc}
The matrix $LT_n^m(a,f)$ is not normal in general, but if we replace in its definition the Toeplitz matrix $T_{\lfloor n/m\rfloor}(f)$ with the circulant matrix $C_{\lfloor n/m\rfloor}(f)$, then we obtain a normal matrix.

\begin{definition}[\textbf{locally circulant operator}]\label{LCop}
Let $m,n\in\mathbb N$ with $m\le n$, let $a:[0,1]\to\mathbb C$ and let $f\in L^1([-\pi,\pi])$. The locally circulant (LC) operator is defined as the following $n\times n$ normal matrix:
\begin{align*}
LC_n^m(a,f)&=D_m(a)\otimes C_{\left\lfloor n/m\right\rfloor}(f)\:\oplus\:O_{n\,{\rm mod}\,m}=\mathop{\rm diag}_{i=1,\ldots,m}\biggl[a\Bigl(\frac im\Bigr)C_{\left\lfloor n/m\right\rfloor}(f)\biggr]\:\oplus\:O_{n\,{\rm mod}\,m}\\
&=\left[\:
\begin{array}{ccccc}
\cline{1-1}
\multicolumn{1}{|c|}{a(\frac1m)C_{\lfloor n/m\rfloor}(f)} &&&& \\
\cline{1-2} 
& \multicolumn{1}{|c|}{a(\frac2m)C_{\lfloor n/m\rfloor}(f)} &&&\\
\cline{2-2}
&& \ddots &&\\
\cline{4-4}
&&&  \multicolumn{1}{|c|}{a(1)C_{\lfloor n/m\rfloor}(f)}
&\\
\cline{4-5}
&&&& \multicolumn{1}{|c|}{O_{n\,{\rm mod}\,m}}
\\\cline{5-5}
\end{array}
\:\right].
\end{align*}
\end{definition}

Note that, by Theorem~\ref{sigmaC},
\begin{equation}\label{LCnm(a,f)}
LC_n^m(a,f)=Q_{n,m}\Delta_{n,m}(a,f^{(\lfloor n/m\rfloor-1)})Q_{n,m}^*,
\end{equation}
where $Q_{n,m}$ is defined in \eqref{Qnm} and $\Delta_{n,m}(a,f)$ is the diagonal matrix
\begin{align*}
\Delta_{n,m}(a,f)&=\left[\:
\begin{array}{ccccc}
\cline{1-1}
\multicolumn{1}{|c|}{a(\frac1m)\Delta_{\lfloor n/m\rfloor}(f^{(\lfloor n/m\rfloor-1)})} &&&& \\
\cline{1-2} 
& \multicolumn{1}{|c|}{a(\frac2m)\Delta_{\lfloor n/m\rfloor}(f^{(\lfloor n/m\rfloor-1)})} &&&\\
\cline{2-2}
&& \ddots &&\\
\cline{4-4}
&&&  \multicolumn{1}{|c|}{a(1)\Delta_{\lfloor n/m\rfloor}(f^{(\lfloor n/m\rfloor-1)})}
&\\
\cline{4-5}
&&&& \multicolumn{1}{|c|}{O_{n\,{\rm mod}\,m}}
\\\cline{5-5}
\end{array}
\:\right].
\end{align*}
For our purposes, we are interested only in $LC_n^m(a,f)$ with $m=m(n)$ such that $m(n)\to\infty$ and $m(n)=o(n)$ as $n\to\infty$.

\begin{lemma}\label{LT=LC}
Let $a:[0,1]\to\mathbb C$ and let $f$ be a trigonometric polynomial. Then,
\[ \{LC_n^{m(n)}(a,f)\}_n\equiv_{\rm a.c.s.}\{LT_n^{m(n)}(a,f)\}_n \]
for every $m(n)\to\infty$ such that $m(n)=o(n)$.
\end{lemma}
\begin{proof}
Throughout this proof, $m=m(n)$ and $m(n)$ has the properties specified in the statement.
Let $r$ be the degree of $f$. By writing explicitly $T_k(f)$ and $C_k(f)$, 
we see that the difference $T_k(f)-C_k(f)$ has at most $2r$ non-zero rows. This implies that ${\rm rank}(T_k(f)-C_k(f))\le2r$ for all $k\in\mathbb N$.
Hence, for every $n$, 
\[ {\rm rank}(LC_n^m(a,f) - LT_n^m(a,f)) \le m\,{\rm rank}(C_{\lfloor n/m\rfloor}(f) - T_{\lfloor n/m\rfloor}(f))\le 2rm = o(n). \]
It follows that $\{LC_n^m(a,f) - LT_n^m(a,f)\}_n$ is zero-distributed and the result is proved.
\end{proof}

Given $m,n\in\mathbb N$ with $m\le n$ and $\zeta:[0,1]\times[0,2\pi]\to\mathbb C$, we denote by $\Delta_{n,m}(\zeta)$ the diagonal matrix
\begin{equation}\label{Deltanmzeta}
\Delta_{n,m}(\zeta)=	\left[\:
	\begin{array}{ccccc}
	\cline{1-1}
	\multicolumn{1}{|c|}{\Delta_{\lfloor n/m\rfloor}(\zeta(\frac1m,\cdot))} &&&& \\
	\cline{1-2} 
	& \multicolumn{1}{|c|}{\Delta_{\lfloor n/m\rfloor}(\zeta(\frac2m,\cdot))} &&&\\
	\cline{2-2}
	&& \ddots &&\\
	\cline{4-4}
	&&&  \multicolumn{1}{|c|}{\Delta_{\lfloor n/m\rfloor}(\zeta(1,\cdot))}
	&\\
	\cline{4-5}
	&&&& \multicolumn{1}{|c|}{O_{n\,{\rm mod}\,m}}
	\\\cline{5-5}
	\end{array}
	\:\right].
\end{equation}

\begin{lemma}\label{spectral_symbol}
Let $\zeta:[0,1]\times[0,2\pi]\to\mathbb C$ be continuous and let $m(n)\to\infty$ such that $m(n)=o(n)$ as $n\to\infty$. Then, $\{\Delta_{n,m(n)}(\zeta)\}_n\sim_\lambda\zeta$.
\end{lemma}
\begin{proof}
Throughout this proof, $m=m(n)$ and $m(n)$ has the properties specified in the statement.
The lemma is proved by direct computation. The eigenvalues of the diagonal matrix $\Delta_{n,m}(\zeta)$ are just the diagonal entries and so they are given by
\[ \zeta\Bigl(\frac im,\frac{2\pi j}{\lfloor n/m\rfloor}\Bigr),\qquad i=1,\ldots,m,\qquad j=1,\ldots,\lfloor n/m\rfloor, \]
plus further $n\,{\rm mod}\,m$ zero eigenvalues. Thus, for every $F\in C_c(\mathbb C)$,
\begin{align*}
\frac1n\sum_{k=1}^nF(\lambda_k(\Delta_{n,m}(\zeta)))&=\frac{(n\,{\rm mod}\,m)F(0)}n+\frac1n\sum_{i=1}^m\sum_{j=1}^{\lfloor n/m\rfloor}F\Bigl(\zeta\Bigl(\frac im,\frac{2\pi j}{\lfloor n/m\rfloor}\Bigr)\Bigr)\\
&=\frac{(n\,{\rm mod}\,m)F(0)}n+\frac{m\lfloor n/m\rfloor}{2\pi n}\cdot\frac{2\pi}{m\lfloor n/m\rfloor}\sum_{i=1}^m\sum_{j=1}^{\lfloor n/m\rfloor}F\Bigl(\zeta\Bigl(\frac im,\frac{2\pi j}{\lfloor n/m\rfloor}\Bigr)\Bigr)\\
&\xrightarrow{n\to\infty}\frac1{2\pi}\int_{[0,1]\times[0,2\pi]}F(\zeta(x,\theta)){\rm d}x{\rm d}\theta,
\end{align*}
where the last limit follows from the assumptions on $m$ and the fact that $\frac{2\pi}{m\lfloor n/m\rfloor}\sum_{i=1}^m\sum_{j=1}^{\lfloor n/m\rfloor}F(\zeta(\frac im,\frac{2\pi j}{\lfloor n/m\rfloor}))$ is a Riemann sum for the continuous (and hence Riemann-integrable) function $F(\zeta(x,\theta))$.
\end{proof}

\subsubsection{Finalization of the proof}

The last results we need to prove Theorem~\ref{normal} are the following two technical lemmas.

\begin{lemma}\label{lem2.8}
Let $\kappa:[0,1]\times[-\pi,\pi]\to\mathbb C$ be measurable. Then, there exists functions $a_{i,\ell}$, $f_{i,\ell}$, $i=1,\ldots,N_\ell$, with $a_{i,\ell}\in C^\infty([0,1])$ and $f_{i,\ell}$ trigonometric polynomial, such that
\[ \sum_{i=1}^{N_\ell}a_{i,\ell}(x)f_{i,\ell}(\theta)\to\kappa(x,\theta)\,\mbox{ a.e.\ in }\,[0,1]\times[-\pi,\pi]. \]
Moreover, if $\kappa$ is real then $a_{i,\ell}$, $f_{i,\ell}$ can be chosen to be real for all $i=1,\ldots,N_\ell$ and all $\ell$.
\end{lemma}
\begin{proof}
See \cite[Lemma~2.8]{GLTbookI}.
\end{proof}

\begin{lemma}\label{mn}
Let $\{A_n\}_n$ be a matrix-sequence and let $\{B_{n,m}\}_n\xrightarrow{\rm a.c.s.}\{A_n\}_n$. Then, there exists a sequence $\{m(n)\}_n$ such that $m(n)\to\infty$ and $\{A_n\}_n\equiv_{\rm a.c.s.}\{B_{n,m(n)}\}_n$.
\end{lemma}
\begin{proof}
Since $\{B_{n,m}\}_n\xrightarrow{\rm a.c.s.}\{A_n\}_n$, for every $m$ there exists $n_m$ such that, for $n\ge n_m$,
\[ A_n=B_{n,m}+R_{n,m}+N_{n,m},\qquad{\rm rank}(R_{n,m})\le c(m)n,\qquad\|N_{n,m}\|\le\omega(m), \]
where $\lim_{m\to\infty}c(m)=\lim_{m\to\infty}\omega(m)=0$.
Without loss of generality, we can assume that the integers $n_m$ have been chosen so that the sequence $\{n_m\}_m$ is strictly increasing. Define
\[ m(n)=\left\{\begin{aligned}
&1, &&\mbox{if }n<n_2,\\
&2, &&\mbox{if }n_2\le n<n_3,\\
&3, &&\mbox{if }n_3\le n<n_4,\\
&\ldots
\end{aligned}\right. \]
It is clear that $m(n)\to\infty$. Moreover, for every $n\ge n_2$ we have $n\ge n_{m(n)}$ and so
\[ A_n=B_{n,m(n)}+R_{n,m(n)}+N_{n,m(n)},\qquad{\rm rank}(R_{n,m(n)})\le c(m(n))n,\qquad\|N_{n,m(n)}\|\le\omega(m(n)). \]
Since $c(m(n))\to0$ and $\omega(m(n))\to0$ (because $m(n)\to\infty$), the matrix-sequence $\{R_{n,m(n)}+N_{n,m(n)}\}_n$ is zero-distributed and we conclude that $\{A_n\}_n\equiv_{\rm a.c.s.}\{B_{n,m(n)}\}_n$.
\end{proof}

\begin{proof}[Proof of Theorem~{\rm\ref{normal}}]
Let $\{A_n\}_n\sim_{\rm GLT}\kappa$, let $\{m(n)\}_n$ be a sequence such that $m(n)\le n$, $m(n)\to\infty$ and $m(n)=o(n)$ as $n\to\infty$, and let $Q_n=Q_{n,m(n)}$, where $Q_{n,m}$ is defined in \eqref{Qnm}. By Lemma~\ref{lem2.8}, there exist functions $a_{i,\ell}$, $f_{i,\ell}$, $i=1,\ldots,N_\ell$, with $a_{i,\ell}\in C^\infty([0,1])$ and $f_{i,\ell}$ trigonometric polynomial, such that $a_{i,\ell}$, $f_{i,\ell}$ are real for all $i,\ell$ if $\kappa$ is real and
\[ \sum_{i=1}^{N_\ell}a_{i,\ell}(x)f_{i,\ell}(\theta)\to\kappa(x,\theta)\,\mbox{ a.e.} \]
By {\bf GLT\,2}\,--\,{\bf GLT\,3}, we have
\[ \left\{\sum_{i=1}^{N_\ell}D_n(a_{i,\ell})T_n(f_{i,\ell})\right\}_n\sim_{\rm GLT}\sum_{i=1}^{N_\ell}a_{i,\ell}(x)f_{i,\ell}(\theta). \]
By {\bf GLT\,6}, we obtain
\[ \left\{\sum_{i=1}^{N_\ell}D_n(a_{i,\ell})T_n(f_{i,\ell})\right\}_n\xrightarrow{\rm a.c.s.}\{A_n\}_n. \]
By Lemmas~\ref{LTacsDT} and~\ref{LT=LC}, for every $i,\ell$ we have
\[ \{D_n(a_{i,\ell})T_n(f_{i,\ell})\}_n\equiv_{\rm a.c.s.}\{LT_n^{m(n)}(a_{i,\ell},f_{i,\ell})\}_n\equiv_{\rm a.c.s.}\{LC_n^{m(n)}(a_{i,\ell},f_{i,\ell})\}_n. \]
Thus,
\[ \left\{\sum_{i=1}^{N_\ell} D_n(a_{i,\ell})T_n(f_{i,\ell})\right\}_n\equiv_{\rm a.c.s.}\left\{\sum_{i=1}^{N_\ell}LC_n^{m(n)}(a_{i,\ell},f_{i,\ell})\right\}_n \]
and, consequently,
\[ \left\{\sum_{i=1}^{N_\ell}LC_n^{m(n)}(a_{i,\ell},f_{i,\ell})\right\}_n\xrightarrow{\rm a.c.s.}\{A_n\}_n. \]
By Lemma~\ref{mn}, there exists a sequence $\{\ell(n)\}_n$ such that $\ell(n)\to\infty$ and 
\[ \{A_n\}_n\equiv_{\rm a.c.s.}\left\{\sum_{i=1}^{N_{\ell(n)}}LC_n^{m(n)}(a_{i,\ell(n)},f_{i,\ell(n)})\right\}_n. \]
Hence,
\begin{equation}\label{eccoci}
\left\{\sum_{i=1}^{N_\ell} LC_n^{m(n)}(a_{i,\ell},f_{i,\ell})\right\}_n\xrightarrow{\rm a.c.s.}\left\{\sum_{i=1}^{N_{\ell(n)}}LC_n^{m(n)}(a_{i,\ell(n)},f_{i,\ell(n)})\right\}_n\equiv_{\rm a.c.s.}\{A_n\}_n.
\end{equation}
By \eqref{LCnm(a,f)},
\begin{equation}\label{eccoci'}
\sum_{i=1}^{N_\ell} LC_n^{m(n)}(a_{i,\ell},f_{i,\ell}) = \sum_{i=1}^{N_\ell} Q_n\Delta_{n,m(n)}(a_{i,\ell},f_{i,\ell})Q_n^*=Q_n\left[\sum_{i=1}^{N_\ell}\Delta_{n,m(n)}(a_{i,\ell},f_{i,\ell})\right]Q_n^*.
\end{equation}
The matrix inside the square brackets is diagonal and we denote it by $D_{n,\ell}$. We also set $D_n=D_{n,\ell(n)}$. By \eqref{eccoci}--\eqref{eccoci'},
\[ \left\{Q_nD_nQ_n^*\right\}_n\equiv_{\rm a.c.s.}\{A_n\}_n. \]
Note that $\{Q_n\}_n$ is a unitary sequence that does not depend on $\{A_n\}_n$, so it only remains to prove that $\{D_n\}_n\sim_\lambda\kappa$.
By definition, the a.c.s.\ convergence is not affected by a unitary base change. Hence, from \eqref{eccoci}--\eqref{eccoci'} we obtain
\begin{equation}\label{.2}
\left\{\sum_{i=1}^{N_\ell}\Delta_{n,m(n)}(a_{i,\ell},f_{i,\ell})\right\}_n\xrightarrow{\rm a.c.s.}\left\{D_n\right\}_n.
\end{equation}
Note that $\sum_{i=1}^{N_\ell}\Delta_{n,m(n)}(a_{i,\ell},f_{i,\ell})=\Delta_{n,m(n)}(\zeta)$, where $\zeta(x,\theta)=\sum_{i=1}^{N_\ell}a_{i,\ell}(x)f_{i,\ell}(\theta)$ is continuous on $[0,1]\times\mathbb R$. 
Thus, Lemma~\ref{spectral_symbol} yields
\begin{equation}\label{.3}
\left\{\sum_{i=1}^{N_\ell}\Delta_{n,m(n)}(a_{i,\ell},f_{i,\ell})\right\}_n\sim_\lambda\sum_{i=1}^{N_\ell}a_{i,\ell}(x)f_{i,\ell}(\theta).
\end{equation}
Note that, according to Lemma~\ref{spectral_symbol}, \eqref{.3} holds if we consider for the function $\sum_{i=1}^{N_\ell}a_{i,\ell}(x)f_{i,\ell}(\theta)$ the domain $[0,1]\times[0,2\pi]$, but actually \eqref{.3} holds even if we replace $[0,1]\times[0,2\pi]$ with $[0,1]\times[-\pi,\pi]$, because $\sum_{i=1}^{N_\ell}a_{i,\ell}(x)f_{i,\ell}(\theta)$ is $2\pi$-periodic in the variable $\theta$.
The thesis $\{D_n\}_n\sim_\lambda\kappa$ now follows from Theorem~\ref{diag2.5}, whose hypotheses are satisfied in view of \eqref{.2}--\eqref{.3} and the fact that $\sum_{i=1}^{N_\ell}a_{i,\ell}(x)f_{i,\ell}(\theta)\to\kappa(x,\theta)$ a.e.
\end{proof}


\subsection{Proof of Theorem~\ref{perturbed-normals}}
The proof of Theorem~\ref{perturbed-normals} requires two preliminary results.
The first one is Theorem~\ref{D=GLT}, which was proved in \cite{spectral_measures}.
The second one is Theorem~\ref{N-cI}, which is proved below and is of interest also in itself.

\subsubsection{GLT sequences and diagonal sequences enjoying a spectral distribution}

The next theorem is the main result of \cite{spectral_measures}.
It plays a central role in the proof of Theorem~\ref{perturbed-normals}.

\begin{theorem}\label{D=GLT}
Let $g:[0,1]\to\mathbb C$ be measurable and let $\{D_n\}_n$ be a diagonal sequence such that $\{D_n\}_n\sim_\lambda g$. Then, there exists a sequence of permutation matrices $\{P_n\}_n$ such that $\{P_nD_nP_n^T\}_n\sim_{\rm GLT}\kappa(x,\theta)=g(x)$.
\end{theorem}
\begin{proof}
See \cite[Theorem~2]{spectral_measures}.
\end{proof}

\subsubsection{Normal sequences: from singular value distribution to spectral distribution}

In what follows, a radial function $G:\mathbb C\to\mathbb C$ is a function of the form $G(z)=H(|z-z_0|)$ with $H:\mathbb R\to\mathbb C$ and $z_0\in\mathbb C$.
It was proved in \cite{radial} that any $F\in C_c(\mathbb C)$ can be arbitrarily approximated in $\infty$-norm by linear combinations of Gaussian functions, which are special cases of radial functions. It follows that the vector space generated by the set of radial functions
\begin{equation*}
\mathscr R=\{H(|z-z_0|):\,H\in C_c(\mathbb R),\ z_0\in\mathbb C\}
\end{equation*}
is dense in $C_c(\mathbb C)$ with respect to the $\infty$-norm. In other words, if $F\in C_c(\mathbb C)$ and $\varepsilon>0$, then there exists a linear combination $L$ of radial functions belonging to the set $\mathscr R$ such that $\|F-L\|_\infty<\varepsilon$.

\begin{theorem}\label{N-cI}
Let $\{A_n\}_n$ be a normal sequence such that $\{A_n-cI_n\}_n\sim_\sigma f-c$ for all $c\in\mathbb C$ and for some measurable $f:D\subset\mathbb R^k\to\mathbb C$ with $0<\mu_k(D)<\infty$. Then $\{A_n\}_n\sim_\lambda f$.
\end{theorem}
\begin{proof}
It is enough to prove that
\begin{equation*}
\lim_{n\to\infty}\frac1n\sum_{i=1}^nG(\lambda_i(A_n))=\frac1{\mu_k(D)}\int_DG(f(\xx)){\rm d}\xx
\end{equation*}
for all radial functions $G$ belonging to the set $\mathscr R$. 
Let $G(z)=H(|z-z_0|)$ be a radial function belonging to the set $\mathscr R$. The matrices $A_n-z_0I_n$ are normal by assumption, and so $\sigma_i(A_n-z_0I_n)=|\lambda_i(A_n-z_0I_n)|$ for all $i=1,\ldots,n$. Using the hypothesis $\{A_n-z_0I_n\}_n\sim_\sigma f-z_0$, we obtain
\begin{align*}
\frac1n\sum_{i=1}^nG(\lambda_i(A_n))&=\frac1n\sum_{i=1}^nH(|\lambda_i(A_n)-z_0|)=\frac1n\sum_{i=1}^nH(|\lambda_i(A_n-z_0I_n)|)\\
&\xrightarrow{n\to\infty}\frac{1}{\mu_k(D)}\int_DH(|f(\xx)-z_0|){\rm d}\xx=\frac{1}{\mu_k(D)}\int_DG(f(\xx)){\rm d}\xx,
\end{align*}
and the theorem is proved.
\end{proof}

\subsubsection{Finalization of the proof}
We have now collected all the necessary ingredients to prove Theorem~\ref{perturbed-normals}.

\begin{proof}[Proof of Theorem~{\rm\ref{perturbed-normals}}]
In the case where $Y_n$ satisfies the first or the second condition in the statement of the theorem, the thesis $\{X_n+Y_n\}_n\sim_\lambda f$ is a direct consequence of \cite[Problem~VI.8.2]{Bhatia}. We prove the thesis in the case where $Y_n$ satisfies the third condition. Let $A_n=X_n+Y_n$ and let $g$ be a rearranged version of $f$ on $[0,1]$, i.e., a measurable function $g:[0,1]\to\mathbb C$ such that
\[ \int_0^1F(g(x)){\rm d}x=\int_0^1F(f(\xx)){\rm d}\xx,\qquad\forall\,F\in C_c(\mathbb C). \]
We remark that such a function $g$ exists by \cite[Lemma~6]{spectral_measures}. We also note that $\{X_n\}_n\sim_\lambda f$ is equivalent to $\{X_n\}_n\sim_\lambda g$.
Since $\{X_n\}_n$ is a normal sequence with $\{X_n\}_n\sim_\lambda g$, there exists a unitary sequence $\{U_n\}_n$ such that $X_n=U_nD_nU_n^*$, where $\{D_n\}_n$ is a diagonal sequence with $\{D_n\}_n\sim_\lambda g$. In view of Theorem~\ref{D=GLT}, by replacing $D_n$ with $P_nD_nP_n^T$ for a suitable permutation matrix $P_n$ (if necessary), we can assume that $\{D_n\}_n\sim_{\rm GLT}g(x)$. Thus,
\[ U_n^*A_nU_n=U_n^*X_nU_n+U_n^*Y_nU_n=D_n+Z_n, \]
where $\{Z_n\}_n=\{U_n^*Y_nU_n\}_n$ is zero-distributed like $\{Y_n\}_n$, and so $\{U_n^*A_nU_n\}_n\sim_{\rm GLT}g(x)$ by {\bf GLT\,2}\,--\,{\bf GLT\,3}. It follows from {\bf GLT\,2}\,--\,{\bf GLT\,3} and the relation $\{I_n\}_n=\{T_n(1)\}_n\sim_{\rm GLT}1$ that $\{U_n^*A_nU_n-cI_n\}_n\sim_{\rm GLT}g(x)-c$ for all $c\in\mathbb C$. We conclude that $\{U_n^*A_nU_n-cI_n\}_n\sim_\sigma g-c$ for all $c\in\mathbb C$ (by {\bf GLT\,1}) and $\{U_n^*A_nU_n\}_n\sim_\lambda g$ (by Theorem~\ref{N-cI}). The latter is equivalent to $\{A_n\}_n\sim_\lambda g$ by definition of spectral distribution. This means that $\{A_n\}_n\sim_\lambda f$ because $g$ is just a rearranged version of $f$.
\end{proof}

\subsection{Proof of Theorem~\ref{GLT1-normal}}

Theorem~\ref{GLT1-normal} is a direct consequence of Theorems~\ref{normal} and~\ref{perturbed-normals}.

\begin{proof}[Proof of Theorem~{\rm\ref{GLT1-normal}}]
Let $\{A_n\}_n\sim_{\rm GLT}\kappa$ and assume that the matrices $A_n$ are normal.
By Theorem~\ref{normal}, there exist a unitary sequence $\{Q_n\}_n$ and a zero-distributed sequence $\{Z_n\}_n$ such that $A_n=Q_nD_nQ_n^*+Z_n$, where $\{D_n\}_n$ is a diagonal sequence and $\{D_n\}_n\sim_\lambda\kappa$. The thesis $\{A_n\}_n\sim_\lambda\kappa$ now follows from Theorem~\ref{perturbed-normals} applied with $X_n=Q_nD_nQ_n^*$ and $Y_n=Z_n$.
\end{proof}

\subsection{Proof of Theorem~\ref{f(GLT)=GLT}}
The proof of Theorem~\ref{f(GLT)=GLT} requires some preliminary results on matrix functions and matrix polynomials as well as on the so-called sparsely unbounded matrix-sequences.


\subsubsection{Matrix functions and matrix polynomials}

We first recall the notion of matrix functions in the case of a diagonalizable matrix.
For more details on matrix functions, see \cite{Higham}.

\begin{definition}[\textbf{function of a diagonalizable matrix}]\label{f(A)}
Let $A\in\mathbb C^{n\times n}$ be a diagonalizable matrix and let $f:\Lambda(A)\to\mathbb C$ be a function defined on the spectrum of $A$, denoted by $\Lambda(A)$. 
We define $f(A)$ as the unique matrix in $\mathbb C^{n\times n}$ such that $f(A)\mathbf u=f(\lambda)\mathbf u$ whenever $A\mathbf u=\lambda\mathbf u$. 
\end{definition}
By Definition~\ref{f(A)}, if 
\[ A=XDX^{-1},\qquad D={\rm diag}(\lambda_1,\ldots,\lambda_n), \]
is a spectral decomposition of the diagonalizable matrix $A$, then a spectral decomposition of $f(A)$ is given by
\[ f(A)=Xf(D)X^{-1},\qquad f(D)={\rm diag}(f(\lambda_1),\ldots,f(\lambda_n)). \]

\begin{definition}[\textbf{complex bivariate polynomial}]
A complex bivariate polynomial is a function $p(x,y):\mathbb C\to\mathbb C$ of the form
\[ p(x,y)=a(x,y)+{\rm i}\,b(x,y), \]
where $a(x,y),b(x,y):\mathbb R^2\to\mathbb R$ are real bivariate polynomials. Note that we write $p(x,y)$ instead of $p(z)$ for convenience ($x$ and $y$ are the real part and the imaginary part of $z$, respectively). Note also that a complex bivariate polynomial $p(x,y)$ should not be confused with a polynomial $q(z)$ of the variable $z$, i.e., a (analytic) function of the form $q(z)=a_0+a_1z+a_2z^2+\ldots+a_mz^m$.
\end{definition}

\begin{definition}[\textbf{complex bivariate polynomial of commuting matrices}]
If $p(x,y):\mathbb C\to\mathbb C$ is a complex bivariate polynomial and $A,B$ are commuting matrices in $\mathbb C^{n\times n}$, then we define $p(A,B)$ as the matrix in $\mathbb C^{n\times n}$ obtained by replacing $x$ with $A$ and $y$ with $B$ in the expression $p(x,y)$ (with the usual convention that a constant $\alpha$ independent of $x$ and $y$ in the expression of $p(x,y)$ must be interpreted as $\alpha x^0y^0$ and replaced with $\alpha I_n$).
Note that this definition is well-posed because $A,B$ commute and so the resulting matrix $p(A,B)$ does not depend on the way in which we write the expression $p(x,y)$.
\end{definition}

For example, if $p(x,y)=2+yx+xy^2+{\rm i}(2xy+yx^3)$ and $A,B$ are commuting matrices in $\mathbb C^{n\times n}$, then $p(A,B)=2I_n+BA+AB^2+{\rm i}(2AB+BA^3)$.
If we write $p(x,y)=2+xy+xy^2+2+{\rm i}(2xy+x^3y)$, the resulting matrix $p(A,B)=2I_n+AB+AB^2+{\rm i}(2AB+A^3B)$ is the same because $A$ and $B$ commute.
In what follows, whenever $f(x,y):\mathbb C\to\mathbb C$ is a complex function written as $f(x,y)$ instead of $f(z)$, it is understood that a writing such as $f(z_0)$, with a unique argument $z_0$, must be interpreted as the complex number obtained by replacing $x$ with ${\rm Re}\,z_0$ and $y$ with ${\rm Im}\,z_0$ in the expression $f(x,y)$.

\begin{lemma}\label{q(A)=q(ReA,ImA)}
Let $A\in\mathbb C^{n\times n}$ be a normal matrix and let $p(x,y):\mathbb C\to\mathbb C$ be a complex bivariate polynomial.
Then, $p(A)$ as given by Definition~{\rm\ref{f(A)}} coincides with $p({\rm Re}\,A,{\rm Im}\,A)$, where
\[ {\rm Re}\,A=\frac{A+A^*}2,\qquad{\rm Im}\,A=\frac{A-A^*}{2{\rm i}} \]
are commuting matrices (because $A$ is normal).
\end{lemma}
\begin{proof}
Since $p(x,y)$ can be written as a linear combination of monic monomials in the variables $x,y$, it is enough to prove the result for such monomials. Let $q(x,y)=x^ry^s$ be a monic monomial in the variables $x,y$.
Consider a spectral decomposition of $A$,
\[ A=QDQ^*,\qquad D={\rm diag}(\lambda_1,\ldots,\lambda_n),\qquad Q\ {\rm unitary}. \]
Note that
\begin{alignat*}{3}
{\rm Re}\,A&=Q({\rm Re}\,D)Q^*, &\qquad{\rm Re}\,D&={\rm diag}({\rm Re}\,\lambda_1,\ldots,{\rm Re}\,\lambda_n),\\
{\rm Im}\,A&=Q({\rm Im}\,D)Q^*, &\qquad{\rm Im}\,D&={\rm diag}({\rm Im}\,\lambda_1,\ldots,{\rm Im}\,\lambda_n).
\end{alignat*}
By Definition~\ref{f(A)}, we have
\begin{align*}
q(A)&=Qq(D)Q^*\\
&=Q\mathop{\rm diag}(q(\lambda_1),\ldots,q(\lambda_n))Q^*\\
&=Q\mathop{\rm diag}(({\rm Re}\,\lambda_1)^r({\rm Im}\,\lambda_1)^s,\ldots,({\rm Re}\,\lambda_n)^r({\rm Im}\,\lambda_n)^s)Q^*\\
&=(Q\mathop{\rm diag}({\rm Re}\,\lambda_1,\ldots,{\rm Re}\,\lambda_n)Q^*)^r(Q\mathop{\rm diag}({\rm Im}\,\lambda_1,\ldots,{\rm Im}\,\lambda_n)Q^*)^s\\
&=({\rm Re}\,A)^r({\rm Im}\,A)^s\\
&=q({\rm Re}\,A,{\rm Im}\,A).\tag*{\qedhere}
\end{align*}
\end{proof}

\subsubsection{Sparsely unbounded matrix-sequences}
We report in this section the notion of sparsely unbounded matrix-sequences and some related results that are necessary for the proof of Theorem~\ref{f(GLT)=GLT}.

\begin{definition}[\textbf{sparsely unbounded matrix-sequence}]\label{s.u.}
A matrix-sequence $\{A_n\}_n$ is said to be sparsely unbounded (s.u.)\ if for every $M>0$ there exists $n_M$ such that, for $n\ge n_M$,
\[ \frac{\#\{i\in\{1,\ldots,n\}:\,\sigma_i(A_n)>M\}}{n}\le r(M), \]
where $\lim_{M\to\infty}r(M)=0$.
\end{definition}

Any matrix-sequence enjoying a singular value distribution as per Definition~\ref{dd} is s.u.

\begin{lemma}\label{p5.4}
If $\{A_n\}_n$ is a matrix-sequence such that $\{A_n\}_n\sim_\sigma f$ for some function $f$ then $\{A_n\}_n$ is s.u.
\end{lemma}
\begin{proof}
See \cite[Proposition~5.4]{GLTbookI}.
\end{proof}

\begin{lemma}\label{su-normal}
Let $\{A_n\}_n$ be an s.u.\ normal sequence. Then, the following property holds: for every $M>0$ there exists $n_M$ such that, for $n\ge n_M$,
\begin{equation}\label{form.s}
A_n=\hat A_{n,M}+\tilde A_{n,M},\qquad{\rm rank}(\hat A_{n,M})\le r(M)n,\qquad \|\tilde A_{n,M}\|\le M,
\end{equation}
where $\lim_{M\to\infty}r(M)=0$, the matrices $\hat A_{n,M}$ and $\tilde A_{n,M}$ are normal, and for all functions $g:\mathbb C\to\mathbb C$ satisfying $g(0)=0$ we have
\[ g(\hat A_{n,M}+\tilde A_{n,M})=g(\hat A_{n,M})+g(\tilde A_{n,M}). \]
\end{lemma}
\begin{proof}
Since the matrices $A_n$ are normal, the singular values $\sigma_i(A_n)$, $i=1,\ldots,n$, coincide with the moduli of the eigenvalues $|\lambda_i(A_n)|$, $i=1,\ldots,n$. Since $\{A_n\}_n$ is s.u., for every $M>0$ there exists $n_M$ such that, for $n\ge n_M$,
\[ \frac{\#\{i\in\{1,\ldots,n\}:\,|\lambda_i(A_n)|>M\}}{n}\le r(M), \]
where $\lim_{M\to\infty}r(M)=0$. Let $A_n=U_n\Lambda_n U_n^*$ be a spectral decomposition of $A_n$. Let $\hat\Lambda_{n,M}$ be the matrix obtained from $\Lambda_n$ by setting to 0 all the eigenvalues of $A_n$ whose modulus is less than or equal to $M$, and let $\tilde\Lambda_{n,M}=\Lambda_n-\hat\Lambda_{n,M}$ be the matrix obtained from $\Lambda_n$ by setting to 0 all the eigenvalues of $A_n$ whose modulus is greater than $M$. Then, for $M>0$ and $n\ge n_M$,
\[ A_n=U_n\Lambda_n U_n^*=U_n\hat\Lambda_{n,M}U_n^*+U_n\tilde\Lambda_{n,M}U_n^*=\hat A_{n,M}+\tilde A_{n,M}, \]
where $\hat A_{n,M}=U_n\hat\Lambda_{n,M}U_n^*$ and $\tilde A_{n,M}=U_n\tilde\Lambda_{n,M}U_n^*$. The matrices $\hat A_{n,M}$, $\tilde A_{n,M}$ constructed in this way are normal, satisfy the properties \eqref{form.s} and, moreover, for all functions $g:\mathbb C\to\mathbb C$ satisfying $g(0)=0$ we have
\[ g(\hat A_{n,M}+\tilde A_{n,M})=g(\hat A_{n,M})+g(\tilde A_{n,M}). \tag*{\qedhere} \]
\end{proof}

\subsubsection{Finalization of the proof}
We have now collected all the necessary ingredients to prove Theorem~\ref{f(GLT)=GLT}.

\begin{proof}[Proof of Theorem~{\rm\ref{f(GLT)=GLT}}]
Let $\{A_n\}_n\sim_{\rm GLT}\kappa$ with $A_n$ normal for every $n$, and let $f(x,y):\mathbb C\to\mathbb C$ be continuous.
For each $M>0$, let $p_{m,M}(x,y):\mathbb C\to\mathbb C$ be a sequence of complex bivariate polynomials that converges uniformly to $f(x,y)$ over $[-M,M]^2$:
\[ \lim_{m\to\infty}\|f-p_{m,M}\|_{\infty,[-M,M]^2}=0. \]
Note that such a sequence exists by the Stone--Weierstrass theorem; see, e.g., \cite[Theorem~7.33]{Rudinino}.
By replacing $p_{m,M}$ with $p_{m,M}+f(0)-p_{m,M}(0)$ if necessary, we can assume, without loss of generality, that $p_{m,M}(0)=f(0)$.
Since any GLT sequence is s.u.\ (by {\bf GLT\,1} and Lemma~\ref{p5.4}), the sequence $\{A_n\}_n$ is s.u.
Hence, by Lemma~\ref{su-normal}, for all $M>0$ there exists $n_M$ such that, for $n\ge n_M$,
\begin{equation}\label{ex8.10}
A_n=\hat A_{n,M}+\tilde A_{n,M},\qquad{\rm rank}(\hat A_{n,M})\le r(M)n,\qquad\|\tilde A_{n,M}\|\le M,
\end{equation}
where $r(M)\to0$ as $M\to\infty$, the matrices $\hat A_{n,M}$ and $\tilde A_{n,M}$ are normal, and for all functions $g:\mathbb C\to\mathbb C$ satisfying $g(0)=0$ we have
\[ g(\hat A_{n,M}+\tilde A_{n,M})=g(\hat A_{n,M})+g(\tilde A_{n,M}). \]
Taking into account that $(f-p_{m,M})(0)=0$, for every $M>0$, every $m$ and every $n\ge n_M$, we can write
\begin{align}\label{almost.done}
f(A_n)&=p_{m,M}(A_n)+f(A_n)-p_{m,M}(A_n)\notag\\
&=p_{m,M}(A_n)+(f-p_{m,M})(\hat A_{n,M})+(f-p_{m,M})(\tilde A_{n,M})\notag\\
&=p_{m,M}(A_n)+R_{n,m,M}+N_{n,m,M},
\end{align}
where, in view of \eqref{ex8.10}, the matrices $R_{n,m,M}=(f-p_{m,M})(\hat A_{n,M})$ and $N_{n,m,M}=(f-p_{m,M})(\tilde A_{n,M})$ satisfy
\begin{align*}
{\rm rank}(R_{n,m,M})&\le{\rm rank}(\hat A_{n,M})\le r(M)n,\\
\|N_{n,m,M}\|&\le\|f-p_{m,M}\|_{\infty,[-M,M]^2}.
\end{align*}
Choose a sequence $\{M_m\}_m$ such that
\begin{equation}\label{crucial}
M_m\to\infty,\qquad\|f-p_{m,M_m}\|_{\infty,[-M_m,M_m]^2}\to0.
\end{equation}
Then, for every $m$ and every $n\ge n_{M_m}$,
\begin{align*}
f(A_n)&=p_{m,M_m}(A_n)+R_{n,m,M_m}+N_{n,m,M_m},\\
{\rm rank}(R_{n,m,M_m})&\le r(M_m)n,\\
\|N_{n,m,M_m}\|&\le\|f-p_{m,M_m}\|_{\infty,[-M_m,M_m]^2},
\end{align*}
which implies that
\[ \{p_{m,M_m}(A_n)\}_n\xrightarrow{\rm a.c.s.}\{f(A_n)\}_n. \] 
Moreover, by Lemma~\ref{q(A)=q(ReA,ImA)}, $p_{m,M_m}(A_n)=p_{m,M_m}({\rm Re}\,A_n,{\rm Im}\,A_n)$, where $\{{\rm Re}\,A_n\}_n\sim_{\rm GLT}{\rm Re}\,\kappa$ and $\{{\rm Im}\,A_n\}_n\sim_{\rm GLT}{\rm Im}\,\kappa$ by {\bf GLT\,3}. Hence, again by {\bf GLT\,3},
\[ \{p_{m,M_m}(A_n)\}_n=\{p_{m,M_m}({\rm Re}\,A_n,{\rm Im}\,A_n)\}_n\sim_{\rm GLT}p_{m,M_m}({\rm Re}\,\kappa,{\rm Im}\,\kappa)=p_{m,M_m}(\kappa). \]
Finally, by \eqref{crucial},
\[ p_{m,M_m}(\kappa)\to f(\kappa)\ \,{\rm a.e.} \]
We conclude that $\{f(A_n)\}_n\sim_{\rm GLT}f(\kappa)$ by {\bf GLT\,5}.
\end{proof}

%


\section*{Acknowledgements}
The authors are members of the Research Group GNCS (Gruppo Nazionale per il Calcolo Scientifico) of INdAM (Istituto Nazionale di Alta Matematica).
Giovanni Barbarino was supported by an Academy of Finland grant (Suomen Akatemian P\"a\"at\"os 331240), by the Alfred Kordelinin S\"a\"ati\"o Grant 210122, and by the ERC Consolidator Grant 101085607 through the Project eLinoR.
Carlo Garoni was supported by the MUR Excellence Department Projects Math@TOV (CUP E83C18000100006) and MatMod@TOV (CUP E83C23000330006) awarded to the Department of Mathematics of the University of Rome Tor Vergata, and by the Department of Mathematics of the University of Rome Tor Vergata through the Project RICH\hspace{1pt}\rule{5pt}{0.4pt}GLT (CUP E83C22001650005).


\begin{thebibliography}{99}

\footnotesize

\bibitem{d-complete}
{\sc Barbarino G.}
{\em Equivalence between GLT sequences and measurable functions.}
Linear Algebra Appl. 529 (2017) 397--412.



\bibitem{spectral_measures}
{\sc Barbarino G.}
{\em Spectral measures.}
Springer INdAM Series 30 (2019) 1--24.

\bibitem{barbarinoREDUCED}
{\sc Barbarino G.}
{\em A systematic approach to reduced GLT.}
BIT Numer. Math. 62 (2022) 681--743.

\bibitem{sdau}
{\sc Barbarino G., Ekstr\"om S.-E., Garoni C., Meadon D., Serra-Capizzano S., Vassalos P.}
{\em From asymptotic distribution and vague convergence to uniform convergence, with numerical applications.}
ArXiv:2309.03662 (2023).

\bibitem{ela}
{\sc Barbarino G., Garoni C.}
{\em From convergence in measure to convergence of matrix-sequences through concave functions and singular values.}
Electron. J. Linear Algebra 32 (2017) 500--513.

\bibitem{rGLT}
{\sc Barbarino G., Garoni C., Mazza M., Serra-Capizzano S.}
{\em Rectangular GLT sequences.}
Electron. Trans. Numer. Anal. 55 (2022) 585--617.

\bibitem{GLTbookIII}
{\sc Barbarino G., Garoni C., Serra-Capizzano S.}
{\em Block generalized locally Toeplitz sequences: theory and applications in the unidimensional case.}
Electron. Trans. Numer. Anal. 53 (2020) 28--112.

\bibitem{GLTbookIV}
{\sc Barbarino G., Garoni C., Serra-Capizzano S.}
{\em Block generalized locally Toeplitz sequences: theory and applications in the multidimensional case.}
Electron. Trans. Numer. Anal. 53 (2020) 113--216.

\bibitem{Bhatia}
{\sc Bhatia R.}
{\em Matrix Analysis.}
Springer, New York (1997).

\bibitem{DavideCalcolo2021}
{\sc Bianchi D.}
{\em Analysis of the spectral symbol associated to discretization schemes of linear self-adjoint differential operators.}
Calcolo 58 (2021) 38.

\bibitem{Ftest}
{\sc Bianchi D., Garoni C.}
{\em On the asymptotic spectral distribution of increasing size matrices: test functions, spectral clustering, and asymptotic estimates of outliers.}
Linear Algebra Appl. (in press).

\bibitem{maximum_norm}
{\sc Bogoya J. M., B\"ottcher A., Grudsky S. M., Maximenko E. A.}
{\em Maximum norm versions of the Szeg\H o and Avram--Parter theorems for Toeplitz matrices.}
J. Approx. Theory 196 (2015) 79--100.

\bibitem{Bottcher-ded-Grudsky}
{\sc Bogoya J. M., B\"ottcher A., Maximenko E. A.}
{\em From convergence in distribution to uniform convergence.}
Bol. Soc. Mat. Mex. 22 (2016) 695--710.

\bibitem{SbMath}
{\sc B\"ottcher A., Garoni C., Serra-Capizzano S.}
{\em Exploration of Toeplitz-like matrices with unbounded symbols is not a purely academic journey.}
Sb. Math. 208 (2017) 1602--1627.

\bibitem{BoSi}
{\sc B\"ottcher A., Silbermann B.}
{\em Introduction to Large Truncated Toeplitz Matrices.}
Springer-Verlag, New York, 1999.

\bibitem{radial}
{\sc Cao F. L., Xie T. F.}
{\em The rate of approximation of Gaussian radial basis neural networks in continuous function space.}
Acta Math. Sin. Engl. Ser. 29 (2013) 295--302.

\bibitem{Davis}
{\sc Davis P. J.}
{\em Circulant Matrices.}
Second Edition, AMS Chelsea Publishing, 1994.

\bibitem{GLTbookI}
{\sc Garoni C., Serra-Capizzano S.}
{\em Generalized Locally Toeplitz Sequences: Theory and Applications (Volume I).}
Springer, Cham, 2017.

\bibitem{GLTbookII}
{\sc Garoni C., Serra-Capizzano S.}
{\em Generalized Locally Toeplitz Sequences: Theory and Applications (Volume II).}
Springer, Cham, 2018.

\bibitem{Higham}
{\sc Higham N. J.}
{\em Functions of Matrices: Theory and Computation.}
SIAM, Philadelphia, 2008.

\bibitem{Rudinino}
{\sc Rudin W.}
{\em Principles of Mathematical Analysis.}
Third Edition, McGraw-Hill, New York (1976).

\end{thebibliography}
\end{document}